\documentclass[12pt]{amsart}
\usepackage{amsthm}
\usepackage[cyr]{aeguill}
\usepackage{dsfont}
\usepackage{fancyhdr}
\usepackage{upgreek}
\usepackage{mathrsfs}
\usepackage{eufrak}
\usepackage{amsmath,amsfonts,amssymb,color,multicol,subeqnarray}\usepackage{graphicx,longtable}
\usepackage[T1]{fontenc}
\usepackage[latin1]{inputenc}  
\usepackage[a4paper]{geometry}
\usepackage[english]{babel}
\theoremstyle{plain}
\newtheorem{teo}{Theorem}
\newtheorem{pro}{Proposition}
\newtheorem{lem}{Lemma}

\theoremstyle{definition}

\theoremstyle{remark}
\newtheorem*{nota}{\bf Remark}

\def\sca #1_#2_#3 {\langle#1,#2\rangle_#3}

\def\C{\mathbb {C}}
\def\m{\mu}
\def\<{\langle}
\def\>{\rangle}

\def\n{\nu}

\def\vp{{\varphi_{\omega}}}
\def\pp{\varphi}

\def\e{\varepsilon}

\font\dsrom=dsrom10 scaled 1200
\def \indic{\textrm{\dsrom{1}}}
\newcommand{\bq}{\begin{equation*}}

\newcommand{\eq}{\end{equation*}}
\newcommand{\ba}{\begin{eqnarray*}}
\newcommand{\ea}{\end{eqnarray*}}

\newcommand{\tr}{\hbox{tr}}

\newcommand{\ove}[1]{\overline{#1}}

\newcommand{\ban}{\begin{eqnarray}}
\newcommand{\ean}{\end{eqnarray}}

\begin{document}

\title{A Lévy-Khinchin Formula for the Space of Infinite Dimensional Square Complex Matrices}
\bigskip
\author{\bf Marouane Rabaoui}


\date{\today} 

\maketitle


\begin{abstract} Using a generalised Bochner type representation for Olshanski spherical pairs, we establish a Lévy-Khinchin formula for the continuous functions of negative type on the space $V_\infty= M(\infty, \mathbb C)$ of infinite dimensional square complex matrices relatively to the action of the product group $K_\infty=U(\infty)\times U(\infty)$. The space $V_\infty$ is the inductive limit of the spaces $V_n=M(n, \mathbb C)$, and the group $K_\infty$ is the inductive limit of the product groups $K_n=U(n)\times U(n)$, where $U(n)$ is the unitary group.\\ 
\end{abstract}

\section{Introduction}

The remarkable progress in convexity theory during the years 1955 to 1965 and thereafter had enabled the progress in the study of infinitely divisible probability measures and the central limit problem. In the pioneering work of G. Choquet (cf. \cite{Cho1}, \cite{Cho2} and \cite{Phel}), the well-known Krein-Milman theorem was extended to an integral representation of elements of a compact convex subset of a locally convex vector space by measures which under additional assumptions are supported by the extreme points of that set. This integral representation theorem yields a new approach to a wealth of problems in analysis such as the Riesz representation, Bernstein's theorem, invariant and ergodic measures, Bochner's theorem of harmonic analysis and the Lévy-Khinchin decomposition of infinitely divisible probability measures. The application of the Choquet theory to the latter two problems proved to be especially useful.\\

The importance of continuous functions of negative type and the Lévy-Khinchin decomposition rests in their applications in the theory of limit theorems for independent and identically distributed random variables, cf. the monograph by B.W. Gnedenko and A.N. Kolmogorov \cite{Gen} or the more recent book by V.V. Petrov \cite{Pet}. The Lévy-Khinchin decomposition can also be used to decompose a given Lévy process into simpler processes (cf. \cite{Boc}). \\

The Lévy-Khinchin decomposition, known also as the integral representation of functions of negative type, or the L{\'e}vy-Khinchin formula has been established in terms of Fourier transforms of measures and functions of negative type for various classes of topological groups. We only mention the references C. Berg, J. P. Christensen and P. Ressel. \cite{Ber1}, C. Berg and G. Forst. \cite{Ber2}, K. R. Parthasarathy \cite{Par}, S. Bochner \cite{Boc}, K. Harzallah \cite{Har} and K. Harzallah and J. Faraut \cite{Far5}.\\

In the frame of G. Olshanski's theory for spherical pairs which had been elaborated in 1990 (cf. \cite{Olsh} and \cite{Far2}), many r{e}sults had been obtained. I. J. Schoenberg, for example, had {e}stablished an analogous of the L{\'e}vy-Khinchin formula for continuous functions of n{e}gative type on $\Bbb R^{(\infty)}$, which are invariant by the orthogonal group $O(\infty)$ (cf. \cite {Sho2}). As for the functions of n{e}gative type d{e}fined on the space of infinite dimensional hermitian matrices and the infinite symmetric group, M. Bouali {e}stablished a L{\'e}vy-Khinchin formula (cf. \cite{BBouali} and \cite{Bouali}).\\

In this paper, we consider the spherical pair $(G_\infty, \ K_\infty)$, which is the inductive limit of the sequence of Gelfand pairs $(G_n, \ K_n)$ : $$G_n=K_n\ltimes V_n, \ \ K_n=U(n)\times U(n), \ \ V_n=M(n, \mathbb C), $$ $$G_\infty=K_\infty\ltimes V_\infty, \ \ K_\infty=U(\infty)\times U(\infty).$$ Here $V_\infty=M(\infty,\mathbb C)$ is the space of infinite complex matrices having only a finite number of non-zero entries, and $U(\infty)$ is the group of the infinite unitary matrices $(u_{ij})$ with complex coefficients such that $u_{ij}=\delta_{ij}$ for $i+j$ large enough. The group $G_\infty$ is equipped with the inductive limit topology. The subgroup $K_\infty$ is closed. The homogeneous space $G_\infty/K_\infty$ is the vector space $$V_\infty=\bigcup_{n=1}^\infty V_n.$$ The law of the inductive limit group $G_\infty=K_\infty\ltimes V_\infty$ is given by : $$(u, \ x)(v, \ y)=\big((u_1v_1, \ u_2v_2), \ x+u_1yu_2^*\big), $$ where $$u=(u_1, u_2), \ v=(v_1, v_2)\in K_\infty \ {\rm and} \ x, \ y\in V_\infty.$$

\bigskip

  Our approach is inspired from a proof given by Berg, Christensen and Ressel in \cite{Ber1} and the work of M. Bouali in \cite{BBouali} and \cite{Bouali}. The main ingredient to establish the Lévy-Khinchin formula for the spherical pair $(G_\infty, \ K_\infty)$ is the generalized Bochner Theorem (cf. \cite{Rab}, Theorem 7). We also make use of sph{e}rical functions of positive type relatively to the spherical pair $(G_\infty, \ K_\infty)$ (cf. \cite{Rab2}).\\
	
	
   Let $G$ be a Hausdorff topological group having $e$ as unit, and $K$ a closed subgroup of $G$. A function $\pp:G\longrightarrow \C$ is said to be of positive type if the kernel defined on $G\times G$ by $(g_1,g_2)\longmapsto\pp(g_2^{-1}g_1)$ is of positive type, i.e. for all $g_1, g_2,\ldots, g_n \in G$ and all $c_1, c_2,\dots, c_n \in \C$, $$\sum_{i=1}^n\sum_{j=1}^nc_i\ove{c_j}\pp(g_j^{-1}g_i)\ \geq 0 .$$ Every function $\pp$ of positive type on $G$ is hermitian, i.e. for all $g\in G$, $\ove{\pp(g)}=\pp(g^{-1})$. In addition, the function $\pp$ is bounded : $|\pp(g)|\leq\pp(e)$. Besides, a function $\pp$ defined on $G$ is said to be $K$-biinvariant if it verifies $\pp(k_1gk_2)=\pp(g)$, for all $k_1$, $k_2\in K$ and all $g\in G$.\\
	
	
    Let $\mathfrak{P}$ be the set of $K_\infty$-biinvariant continuous functions of positive type on $G_\infty$ satisfying $\pp(0)=1$. The set ${\rm ext(}\mathfrak{P}{\rm )}$ of extreme points of this convex cone are indeed the spherical functions of positive type relatively to the pair $(G_\infty, \ K_\infty)$. A function $\pp$ on $G_\infty$ which is right invariant under $K_\infty$ does not depend on the variable $u\in K_\infty$. Therefore, it is possible to see it as a function on $V_\infty$ : $\pp(g)=\pp\big((u, x)\big)=\pp_0(x).$ Moreover, if the function $\pp$ is $K_\infty$-biinvariant then the function $\pp_0$ which is defined on $V_\infty$ is $K_\infty$-invariant, or $U(\infty)$-biinvariant. Also, the function $\pp$ is of positive type on $G_\infty$ if and only if the function $\pp_0$ is of positive type. \\
		

The topology defined on ${\rm ext(}\mathfrak{P}{\rm )}$ can be expressed in terms of a set $\Omega$ of parameters. This enables us to prove a parameterized version of the generalized Bochner theorem (\cite{Rab}, Theorem 7), which represents the key theorem used, in the third section of this paper, to get a Lévy-khinchin formula for the space $V_\infty^2$ of all square complex matrices of Hilbert-Schmidt on which the group $K_\infty $ acts: \\


\begin{paragraph}{\bfseries{Main Theorem}} \ \ Let $\psi$ be a continuous $K_{\infty}$-invariant function on $V_\infty^2$. Then $\psi$ is of negative type if and only if it has the following integral repr{e}sentation 
$$\psi(\xi)=\psi(0)+A\;
{\rm tr}(\xi^2)+\int_{\Omega\backslash\{0\}}(1-\varphi_{\omega}(\xi))\,\n(d\omega),$$
where $A$ is a positive constant, $\varphi_{\omega}$ is the spherical function associated to $\omega$ and $\n$ is a positive measure on $\Omega\setminus\{0\}$ such that $$\int_{\Omega\backslash\{0\}}\frac{||\omega||}{1+||\omega||}\;\n(d\omega)<\infty.$$ 

\vspace{-0.05cm}

The constant $A$ and the measure $\n$ are uniquely determin{e}d. \\ \\
\end{paragraph}

\section{Parameterization of spherical functions on $(G_\infty, \ K_\infty)$}


\bigskip

Let $D_\infty$ denote the subspace of diagonal matrices in $V_\infty$. An element of $D_\infty$ can be decomposed as diag$(a_1, a_2, \dots)$, with $a_1,a_2,\dots\in\mathbb R$ and $a_j=0$ for $j$ large enough. Any matrix $x\in V_\infty$ can be diagonalised as $$x=u\ {\rm diag}(a_1, a_2, \dots) \ v^* \ \ \big(u, v\in U(\infty)\big).$$ Consequently, any $K_\infty$-invariant function on $V_\infty$ is uniquely determined by its restriction to the subset $D_\infty$. The spherical functions of positive type relatively to the pair $(G_\infty, \ K_\infty)$ denoted as $\pp_\omega\in{\rm ext(}\mathfrak{P}{\rm )}$ are given for all $\xi\in V_\infty$ by (cf. \cite{Rab2}) : $$\pp_\omega(\xi)=\pp_\omega\big({\rm diag}(\xi_1, \dots, \xi_n, 0, \dots)\big)=\Pi(\omega,\xi_1)\dots\Pi(\omega,\xi_n)={\rm det}\;\Pi(\omega,\xi), $$ where $\Pi(\omega,.)$ is called the modified P{\'o}lya function of parameter $\omega=(\alpha,\gamma)$ and defined on $\mathbb R$ by : $$\Pi(\omega,\lambda):=e^{-\gamma\lambda^2}\prod_{k=1}^\infty\frac{1}{1+\alpha_k\lambda^2}, $$ with $$\omega=(\alpha,\gamma),\ \gamma\in\mathbb R_+, \ \alpha_k\in\mathbb R_+ \ {\rm and} \ \sum_{k=1}^\infty \alpha_k <\infty.$$

We consider on the set $\mathfrak{B}$ of modified P{\'o}lya functions the topology of uniform convergence on compact sets of $\mathbb R$. The topological space $\mathfrak{B}$ is metrizable and complete. This topology can be expressed in terms of the set of parameters : $$\Omega=\left\{\omega=(\alpha,\gamma)\,\Big|\, \alpha=(\alpha_j)_{j\geq1},\, \alpha_j\geq0,\, \sum_{j=1}^\infty\alpha_j<\infty,\, \gamma\geq0\right\}.$$ For a continuous function $f$ on $\mathbb R$, we define the function $L_f$ on $\Omega$ by \bq\label{log}L_f(\omega)=\int_{\mathbb R} f(t)\sigma_\omega(dt)=\gamma f(0)+\sum\limits_{j=1}^{\infty}
{\alpha}_j f(\alpha_j).\eq Let us remark that the moments of the measure $\sigma_\omega$ are given by $${\mathcal M}_0(\sigma_\omega)=\int_{\mathbb
  R}\sigma_\omega(dt)=\gamma+\sum_{k=1}^{\infty}\alpha_k=\gamma+p_1(\alpha),$$
and for $m\geq 1$,
$${\mathcal M}_m(\sigma_\omega)=\int_{\mathbb
  R}t^m\sigma_\omega(dt)=\sum_{k=1}^{\infty}\alpha_k^{m+1}=p_{m+1}(\alpha),$$
where $p_m$ is the Newton power sum function : for $x=(x_1,x_2,\ldots)\in\ell^{1}(\mathbb N)$ and $m\geq 1$, $$p_m(x)=\sum_{k=1}^{\infty}x_k^m.$$

\bigskip

We consider on $\Omega$ the initial topology associated to the functions $L_f$. A point $\omega\in \Omega$ is seen as a point configuration, i.e. a permutation of the numbers $\{\alpha_k\}$, $\gamma$ does not change $\omega$.

\begin{pro} The correspondence $\Omega\leftrightarrow$ {\rm ext(}$\mathfrak P${\rm)} is an isomorphism between two standard spaces.\\
\begin{paragraph}{\bfseries{Proof.}}{\rm For $\lambda$ fixed, the function $\omega\mapsto \Pi(\omega,\lambda)$ is injective and continuous on $\Omega$. This can be seen by looking at the logarithmic derivative of $\Pi(\omega,\lambda)$ : \bq\label{logg}\begin{aligned}\frac{\Pi^{'}(\omega,\lambda)}{\Pi(\omega,\lambda)}&=-2\;(\gamma+p_{1}(\alpha))\lambda+\sum\limits_{m=2}^{\infty}p_m(-\alpha)\lambda^{2m-1}.\end{aligned}\eq  Besides, the topology of uniform convergence on compact sets of $\mathbb R$ defined on $\mathfrak{B}$ is equivalent to the topology defined on the set of parameters $\Omega$ (cf. \cite{Rab2}, Proposition 4.5). In consequence, the application $\omega\longmapsto {\Pi}(\omega,.)$ defines a homeomorphism between $\Omega$ and $\mathfrak B$.\\ 

 The set $\Omega$ is separable, metrizable and complete. This can be easily deduced from (\cite{Par}, Chapter II, Theorem 6.2). Let us put for every modified P{\'o}lya function of parameter \ $\omega=(\alpha,\gamma)$, $$||\omega||=p_1(\alpha)+\gamma,$$ and, for every $R\geq0$, $$\Omega_R=\{\omega\in\Omega \ | \ ||\omega||\leq R\}.$$ By (\cite{Rab2}, Corollary 4.6), the set $\Omega_R$ is a compact subset of $\Omega$. Since the application $\omega\longmapsto ||\omega|| $ is continuous, the set \ $\left\{\omega\in\Omega\;\mid||\omega||< R\right\}$ \ is open. Hence, every point $\omega$ has a compact neighborhood which proves that $\Omega$ is locally compact. The set $\Omega$ represents, in consequence, a standard space. In addition, the proof of the generalized Bochner theorem (\cite{Rab}, Theorem 7), enables us to see that ext($\mathfrak P$) is a Borel set in a standard space. Hence, it is standard by (\cite{D2}, Appendice B, B 20). Furthermore, the correspondence ext($\mathfrak P$)$\rightarrow\Omega$, $\pp_\omega\mapsto \omega$ is Borelian and one-to-one. In consequence, by (\cite{D2}, Appendice B, B22), it is an isomorphism between two standard spaces.$\hspace{1cm}\Box$\\}
\end{paragraph}
\end{pro}

Hence, we can get a parameterized version of the generalized Bochner theorem (\cite{Rab}, Theorem 7) :\\

\begin{teo}\label{th3}{\rm (Generalized Bochner)} \ Let $\varphi$ be a $K_\infty$-invariant continuous function of positive type on $V_\infty$ with $\varphi(0)=1$. Then, there exists a unique probability measure $\m$ defined on $\Omega$ such that, for every $\xi\in V_\infty$,
\begin{equation}\label{gim}\varphi(\xi)=\int_{\Omega}\varphi_{\omega}(\xi)\, \m(d\omega).\end{equation}
\end{teo}

\section{Analytical and integral representation properties for functions in ${\rm ext(}\mathfrak{P}{\rm )}$}

\bigskip

Let $V_\infty^2$ be the space of all square complex matrices of Hilbert-Schmidt. Its topology is defined via the norm $$|||\xi|||=\left(\sum_{i,j=1}^\infty|\xi_{ij}|^2\right)^{\frac{1}{2}}.$$

Let us recall that a function $\vp\in{\rm ext(}\mathfrak{P}{\rm )}$ is defined on $V_\infty$ by
$$\varphi_{\omega}(\xi)={\rm det}\;\Pi(\omega,\xi)=e^{-\gamma\tr(\xi^2)}\prod_{k=1}^{\infty}\frac
  {1}{\det(1+\alpha_k\xi^2)},$$ with $$\omega=(\alpha,\gamma),\ \gamma\in\mathbb R_+, \ \alpha_k\in\mathbb R_+ \ {\rm and} \ \sum_{k=1}^\infty \alpha_k <\infty.\\$$

\begin{pro}\label{prolong} The function $\vp$ defined on $V_\infty$ can be extended into a continuous function on $V_\infty^2$.\end{pro}
\begin{paragraph}{\bfseries{Proof.}}{\rm  

Since the function $\xi\longmapsto e^{-\gamma\tr(\xi^2)}$ is continuous on $V_{\infty}^2$, it remains to prove that the infinite product is continuous. Let us recall that if $\xi$ is a Hilbert-Schmidt operator then $\xi^2$ is a trace class operator. In consequence, the function $\xi\longmapsto\det(1+\xi^2)$ is continuous on the space of Hilbert-Schmidt operators. Let $\rho<1$, then for $|||\eta|||\leq\rho<1$ we have $$\det(1-\eta)=\exp\tr\left(\log(1-\eta)\right).$$ In consequence, 

$$\begin{aligned}\det\left(1+\eta^2\right)^{-1}&=&\exp-\left(\tr\left(\sum\limits_{m=1}^{\infty}\frac{(-1)^{m}}{m}\;\eta^{2m}\right)\right)\\&=&\exp\left(\sum\limits_{m=1}^{\infty}\frac{(-1)^{m+1}}{m}\;\tr\left(\eta^{2m}\right)\right).\end{aligned}\\$$ Let us note that for a Hilbert-Schmidt operator $\eta$, $$\begin{aligned}|\tr(\eta^{2m})|\leq|||\eta|||^{2m},\end{aligned}$$
and so, $$\begin{aligned}\left|\sum\limits_{m=1}^{\infty}\frac{(-1)^{m+1}}{m}\;\tr\left(\eta^{2m}\right)\right|\leq\sum\limits_{m=1}^{\infty}\frac{|||\eta|||^{2m}}{m}\leq \frac{|||\eta|||^2}{(1-\rho)}, \end{aligned}$$
It follows that there exists a constant $C=C(\rho)$ such that \begin{equation}\label{gom}\det\left(1+\eta^2\right)^{-1}-1\leq\;C\;|||\eta|||^2  .\end{equation}

Let $R>0$. Since $\sum_{k=1}^\infty \alpha_k <\infty$, there exists $k_0\in\mathbb N^*$ such that for all $k\geq k_0$, we have \ $|\alpha_k |\leq\dfrac{\rho}{R^2}$. Therefore, if $|||\xi|||\leq R$, then \ $|||\alpha_k\xi^2 |||\leq
  \rho.$ It follows by equation (\ref{gom}) that  $$\left|\det\left(1+\alpha_k\xi^2\right)^{-1}-1\right|\leq\;C\;\alpha_k\;|||\xi|||^2   .$$

Since $\sum_{k=1}^\infty \alpha_k <\infty$, the infinite product $\prod_{k\geq1}\det(1+\alpha_k\xi^2)^{-1}$ is uniformly convergent in the Hilbert-Schmidt norm on every ball $B(0,R)$ in $V_\infty^2$. Hence, it is a continuous function on $V_\infty^2.\hfill\square\\$

}

\end{paragraph}

\begin{nota}  We can conclude, by using the previous proposition and the dominated convergence theorem, that the function $\varphi$ in formula (\ref{gim}) is also continuous on $V_{\infty}^2$ and so Theorem \ref{th3} holds for $\xi\in V_{\infty}^2$.

\end{nota}

The following results will be necessary in the proof of the Lévy-Khinchin representation. 

\begin{lem}\label{lemm} Any function $\vp\in{\rm ext(}\mathfrak{P}{\rm )}$ can be written for $\omega$ sufficiently close to zero as
$$\vp(\xi)=1-||\omega||{\tr}(\xi^2)+R(\omega,\xi),$$
 where $$\lim\limits_{\omega\to 0}\frac{R(\omega,\xi)}{||\omega||}=0 \ \ (\xi\in V_\infty).\\$$
 
 In addition, for all $\rho>0$, there exist $C>0$ and $\varepsilon >0$, such that, if $||\omega||\leq\varepsilon$ and $|||\xi|||\leq\rho$, then $$|1-\vp(\xi)|\leq C||\omega||.$$
\end{lem}

\begin{paragraph}{\bfseries{Proof.}}{\rm

   ${\rm\bf a)}$ If $\xi=0$, the result is obvious. Let us assume that $\xi\neq 0$. For $\omega\in\Omega$ sufficiently close to zero, we have $$\vp(\xi)=\exp\left(-{\gamma}{\tr}(\xi^2)-\sum\limits_{k=1}^{\infty}{\tr}\log
  {(1+\alpha_k\xi^2)}\right),$$ Hence, for all $\xi\in V_\infty^2$ and all $\omega$ in a neighborhood of zero, we have

\begin{equation}\label{ere}\vp(\xi)=\exp\left(-\gamma {\tr}(\xi^2)-\sum\limits_{k=1}^{\infty}{\tr}\left(\sum_{m=1}^{\infty}\frac{(-1)^{m+1}}{m}\alpha_k^m\xi^{2m}\right)\right).\end{equation}
Let us remark that for all $m\geq 1$
\begin{equation}\label{80}p_m(\alpha)\leq
  ||\omega||^m.\end{equation}

In consequence, for $\omega$ sufficiently small and $||\omega||<\frac{1}{|||\xi|||^2}$, $$\begin{aligned}\sum_{m\geq 1}\sum_{k\geq
  1}\frac{{\alpha_k^m}}{m}|||\xi|||^{2m}&\leq\sum_{m\geq
  1}\left(||\omega||\;|||\xi|||^2\right)^m\; <+\infty.\end{aligned}$$
  
We can then permute the sums in the equation (\ref{ere}) to get
\begin{equation}\label{exp}\vp(\xi)=\exp\left(-||\omega||{\tr}(\xi^2)-\sum_{m=2}^{\infty}\frac{(-1)^{m+1}}{m}p_m(\alpha){\tr}(\xi^{2m})\right).\end{equation}
 Let us put  $$f(\omega,\xi)=-||\omega||{\tr}(\xi^2)-\sum_{m=2}^{\infty}\frac{(-1)^{m+1}}{m}p_m(\alpha){\tr}(\xi^{2m}),$$
and
$$g(\omega,\xi)=-\sum_{m=2}^{\infty}\frac{(-1)^{m+1}}{m}p_m(\alpha){\tr}(\xi^{2m}).$$ 

The equation (\ref{exp}) implies that, \begin{equation}\label{tw}\vp(\xi)=1-||\omega||{\tr}(\xi^2)+g(\omega,\xi)+\sum_{n=2}^{\infty}\frac{(f(\omega,\xi))^n}{n!}.\end{equation}
Using the equation (\ref{80}) and the fact that $|\tr(\xi^{2m})|\leq |||\xi|||^{2m}$, we get
$$|g(\omega,\xi)|\leq\sum_{m=2}^{\infty}(||\omega||\;|||\xi|||^2)^m.$$

It follows that for $\displaystyle||\omega||<\frac{1}{|||\xi|||^2}$ we have
\begin{equation}\label{www}|g(\omega,\xi)|\leq\frac{|||\xi|||^4}{1-||\omega||\;|||\xi|||^2}\;||\omega||^2.\end{equation}
In consequence,
$$|f(\omega,\xi)|\leq\left(|||\xi|||^2+\frac{||\omega||\;|||\xi|||^4}{1-||\omega||\;|||\xi|||^2}\right)||\omega||.$$
Let us put $$ B(\omega,\xi)=|||\xi|||^2+\frac{||\omega||\;|||\xi|||^4}{1-||\omega||\;|||\xi|||^2}.$$
Then,
$$\begin{aligned}\left|\sum_{n=2}^{\infty}\frac{(f(\omega,\xi))^n}{n!}\right|&\leq
\sum_{n=2}^{\infty}\frac{(B(\omega,\xi))^n||\omega||^{n}}{n!}=
||\omega||^2\;\sum_{n=2}^{\infty}\frac{(B(\omega,\xi))^n||\omega||^{n-2}}{n!}.\end{aligned}$$ It follows that if
$||\omega||<\inf\left(\frac{1}{|||\xi|||^2},1\right)$, then
\begin{equation}\label{wwww}\left|\sum_{n=2}^{\infty}\frac{(f(\omega,\xi))^n}{n!}\right|\leq
  ||\omega||^2\;\exp(B(\omega,\xi)).\end{equation}
Let us put
$$R(\omega,\xi)=g(\omega,\xi)+\sum_{n=2}^{\infty}\frac{(f(\omega,\xi))^n}{n!}.$$
By the equations (\ref{www}) and (\ref{wwww}), if $||\omega||<\inf\left(\frac{1}{|||\xi|||^2},1\right)$ and $\omega$ is sufficiently small, we get
\begin{equation}\label{ieiee}|R(\omega,\xi)|\leq
  \left(\exp(B(\omega,\xi))+\frac{|||\xi|||^4}{1-||\omega||\;|||\xi|||^2}\right)||\omega||^2.\end{equation}
Hence, 
 $$\lim\limits_{\omega\to 0}\frac{R(\omega,\xi)}{||\omega||}=0.$$

${\rm\bf b)}$ Let $\rho>0$ and $\displaystyle\e<\inf(\frac{1}{\rho},1)$. If $||\omega||\leq\e$ and
$|||\xi|||\leq\rho$, then
 $$\begin{aligned}
 B(\omega,\xi)=|||\xi|||^2+\frac{||\omega||\;|||\xi|||^4}{1-||\omega||\;|||\xi|||^2}\leq
 \rho^2+\frac{\e \rho^4}{1-\e\rho^2}=C_1.\end{aligned}$$
Using the equation (\ref{ieiee}), we deduce that
$$|R(\omega,\xi)|\leq C_2||\omega||, \quad {\rm where}\quad C_2=\e (\exp(C_1)+C_1-\rho^2).$$
It follows that
$$|1-\vp(\xi)|\leq|||\xi|||^2\;||\omega||+|R(\omega,\xi)|\leq C||\omega||, \quad {\rm where}\quad C=\rho^2+C_2 \ . \qquad\hfill\square\\$$ }

\end{paragraph}

\begin{pro}\label{unique} Let $\nu_1$ et $\nu_2$ be two positive measures on $\Omega\backslash\{0\}$ such that, for all \;$\xi\in V_\infty$,

 \begin{equation}\label{onivom}\int_{\Omega\backslash\{0\}}(1-\vp(\xi))\,\nu_1(d\omega)=\int_{\Omega\backslash\{0\}}(1-\vp(\xi))\,\nu_2(d\omega),\end{equation} and \begin{equation}\label{oon}\int_{\Omega\backslash\{0\}}\frac{||\omega||}{1+||\omega||}\;\nu_i(d\omega)<\infty,\qquad(i=1,2).\end{equation}
 
 \bigskip

Then, for all \;$\xi,\eta\in V_\infty$,

$$\int_{\Omega\backslash\{0\}}\vp(\eta)(1-\vp(\xi))\,\nu_1(d\omega)=\int_{\Omega\backslash\{0\}}\vp(\eta)(1-\vp(\xi))\,\nu_2(d\omega).$$
\end{pro}

\begin{paragraph}{\bfseries{Proof.}}{\rm
The spherical function $\vp$ which is defined on the spherical pair \\ $\displaystyle{(K_{\infty}\ltimes V_\infty,K_\infty)}$ verifies, for all $\xi,\eta\in V_\infty$,
\begin{equation}\label{mni}\lim_{n\rightarrow\infty}\int_{U(n)\times U(n)}\vp(\xi+k_1\eta k_2^*)\,\alpha_n(dk_1)\alpha_n(dk_2)=\vp(\xi)\vp(\eta),\end{equation}
where $\alpha_n$ is the normalized Haar measure of the unitary group $U(n)$.\\

First, we will prove that
\ba\lim_{n\to\infty}\int_{U(n)\times U(n)}\int_{\Omega\backslash\{0\}}(1-\vp(\xi+k_1\eta k_2^*))\;\nu_1(d\omega)\;\alpha_n(dk_1)\alpha_n(dk_2)&\\=\int_{\Omega\backslash\{0\}}(1-\vp(\xi)\vp(\eta))\,\nu_1(d\omega).\qquad\qquad\qquad\qquad\qquad\qquad \ \ \ea

  Let us begin by majorizing the function $$\displaystyle(\omega, k_1, k_2)\mapsto
  (1+||\omega||)\;\frac{1-\vp(\xi+k_1\eta k_2^*)}{||\omega||}.$$
  
Let $R>0$ such that $|||\xi|||+|||\eta|||\leq R$ and $0<\e<\inf(\frac{1}{R},1)$.\\

$a)$ If $||\omega||\leq\sqrt\e$, since $|||\xi+k_1\eta k_2^*|||\leq
|||\xi|||+|||\eta|||\leq R$, then by Lemma \ref{lemm},
$$(1+||\omega||)\frac{|1-\vp(\xi+k_1\eta k_2^*)|}{||\omega||}\leq
(1+\varepsilon)C,$$
where $C$ is a constant d{e}pending on $R$ and $\e$.

$b)$ If $||\omega||\geq\sqrt\e$ then $$(1+||\omega||)\frac{|1-\vp(\xi+k_1\eta k_2^*)|}{||\omega||}\leq 2(1+\frac{1}{\e}).$$

Since the measures $\alpha_n$ and $\displaystyle\frac{||\omega||}{1+||\omega||}\nu_1(d\omega)$ are bounded, the function
$$(\omega, k_1, k_2)\mapsto(1+||\omega||)\frac{|1-\vp(\xi+k_1\eta k_2^*)|}{||\omega||},$$
is int{e}grable with respect to the product measure $$\displaystyle\frac{||\omega||}{1+||\omega||}\nu_1(d\omega)\times
\alpha_n(dk_1)\times\alpha_n(dk_2).$$ Hence, the Fubini theorem implies that

\ba\int_{U(n)\times U(n)}\int_{\Omega\backslash\{0\}}(1-\vp(\xi+k_1\eta
k_2^*))\,\nu_1(d\omega)\,\alpha_n(dk_1)\,\alpha_n(dk_2)\\=\int_{\Omega\backslash\{0\}}\int_{U(n)\times U(n)}(1-\vp(\xi+k_1\eta k_2^*))\,\alpha_n(dk_1)\alpha_n(dk_2)\nu_1(d\omega) .\ea
Using the inequalities in $a)$ and $b)$ besides to the fact that the measure $\alpha_n$ is a probability measure, we conclude that the function

$$\omega\mapsto\int_{U(n)\times U(n)}(1+||\omega||)\frac{1-\vp(\xi+k_1\eta k_2^*)}{||\omega||}\,\alpha_n(dk_1)\,\alpha_n(dk_2),$$ is dominated ind{e}pendently of $n$ and $\omega$. Moreover, the measure $\displaystyle\frac{||\omega||}{1+||\omega||}\nu_1(d\omega)$ is positive and bounded, which implies, by using the dominated convergence theorem and the equation (\ref{mni}), that
\ba\lim_{n\to\infty}\int_{\Omega\backslash\{0\}}\int_{U(n)\times U(n)}(1-\vp(\xi+k_1\eta k_2^*))\,\alpha_n(dk_1)\,\alpha_n(dk_2)\nu_1(d\omega)\\=\int_{\Omega\backslash\{0\}}(1-\vp(\xi)\vp(\eta))\nu_1(d\omega),\qquad\qquad\qquad\qquad\qquad\qquad \ \ \ea
and so,
\ba\lim_{n\to\infty}\int_{U(n)\times U(n)}\int_{\Omega\backslash\{0\}}(1-\vp(\xi+k_1\eta k_2^*))\,\nu_1(d\omega)\,\alpha_n(dk_1)\,\alpha_n(dk_2)\\=\int_{\Omega\backslash\{0\}}(1-\vp(\xi)\vp(\eta))\,\nu_1(d\omega).\qquad\qquad\qquad\qquad\qquad\qquad \ \ \ea

In consequence, the equation (\ref{onivom}) implies that
$$\int_{\Omega\backslash\{0\}}(1-\vp(\xi)\vp(\eta))\,\nu_1(d\omega)=\int_{\Omega\backslash\{0\}}(1-\vp(\xi)\vp(\eta))\,\nu_2(d\omega).$$
The result follows by substituting $\xi$ by $\eta$ in the equation (\ref{onivom}), and considering the diff{e}rence with the previous equation.$\hfill\square$\\}
\end{paragraph}

\section{The L{\'e}vy-Khinchin Formula}


A function $\psi$ d{e}fined on a real vector space $\mathcal V$ with complex values is said to be {\it of negative type} if $\psi(0)\geq0,$ $\psi(-\xi)=\ove{\psi(\xi)}$ and, for all $\xi_1,\dots,\xi_N\in \mathcal V$ and all $c_1,\dots,c_N\in\Bbb C$ such that $\sum_{i=1}^Nc_i=0$, $$\sum_{i,j=1}^Nc_i\ove{c_j}\psi(\xi_i-\xi_j)\leq0.$$  A function of negative type $\psi$ is said to be normalised if $\psi(0)=0$. Since $\psi(\xi)-\psi(0)$ is also continuous, of negative type and biinvariant by $U(\infty)$, we can assume that $\psi(0)=0$. If $\varphi$ is a function of positive type, then $\psi(\xi)=\varphi(0)-\varphi(\xi)$ is of negative type. Moreover, the functions of negative type and those of positive type are related by the following property : \\ 
\begin{pro}\label{pr}{\rm(Schoenberg (\cite{Sho1}, page 527) and (\cite{Ber2}, Theorem 7.8))}\\ The function $\psi$ is of negative type if and only if 
$\psi(0)\geq 0,$ and, for all $t\geq 0$, $e^{-t\psi}$ is of positive type.\\
\end{pro}

The function $\xi\mapsto\pp_\omega(\xi)={\rm det}\;\Pi(\omega, \xi)$, which is d{e}fined on $V_\infty$, can be extended naturally into a continuous function on $V_\infty^2$ (Proposition \ref{prolong}). Using the generalized Bochner theorem (Theorem \ref{th3}), we will establish in the following theorem that every continuous function of negative type on $V_\infty^2$ has an integral representation analogous to the L{\'e}vy-Khinchin formula. The proof uses the same method as in Bouali's paper \cite{BBouali}.\\

\begin{teo}\label{th2}{\rm(L{\'e}vy-Khinchin Formula)}.
Let $\psi$ be a continuous $K_{\infty}$-invariant function on $V_\infty^2$. Then $\psi$ is of negative type if and only if it has the following integral repr{e}sentation 
$$\psi(\xi)=\psi(0)+A\;
{\rm tr}(\xi^2)+\int_{\Omega\backslash\{0\}}(1-\varphi_{\omega}(\xi))\,\n(d\omega),$$
where $A$ is a positive constant and $\n$ is a positive measure on $\Omega\setminus\{0\}$ such that $$\int_{\Omega\backslash\{0\}}\frac{||\omega||}{1+||\omega||}\;\n(d\omega)<\infty.$$ The constant $A$ and the measure $\n$ are uniquely determin{e}d. \end{teo}

\begin{paragraph}{\bfseries{Proof.}}{\rm

${\rm\bf a)}$ \ First, let us remark that, by Lemma \ref{lemm}, the previous integral is well defined and that a function $\psi$ given under the previous representation is invariant by $K_{\infty}$ and of negative type. To prove that it is a continuous function, we use Proposition \ref{prolong} and the dominated convergence theorem.\\

${\rm\bf b)}$ \ {\bf Existence of the representation} : Let $\psi$ be a continuous function of negative type on $V_\infty^2$, which is invariant by $K_\infty$. Since $\psi(\xi)-\psi(0)$ is also $K_\infty$ invariant, continuous and of negative type, then we can assume that $\psi(0)=0$. For $t\geq 0$, the function $e^{-t\psi}$ is $K_\infty$ invariant, continuous and of positive type on $V_\infty^2$. Hence, by the generalized Bochner theorem (Theorem \ref{th3}) and Proposition \ref{prolong}, there exists a unique probability measure $m_t$ on $\Omega$ such that $$e^{-t\psi(\xi)}=\int_{\Omega}\varphi_{\omega}(\xi)\;m_t(d\omega).$$ In consequence, for all $t> 0$, $$\frac{1-e^{-t\psi(\xi)}}{t}=\int_{\Omega}(1-\varphi_{\omega}(\xi))\;\frac{m_t}{t}(d\omega),$$
and we have $$\lim_{t\to
  0}\frac{1-e^{-t\psi(\xi)}}{t}=\psi(\xi)\qquad {\rm and}\qquad\lim_{t\to
  +\infty}\frac{1-e^{-t\psi(\xi)}}{t}=0.$$
In addition, for $\xi$ fixed, there exists a constant $C(\xi)\geq 0$ such that $$0\leq\frac{1-e^{-t\psi(\xi)}}{t}\leq C(\xi),$$
and so
$$\int_{\Omega}(1-\varphi_{\omega}(\xi))\;\frac{m_t}{t}(d\omega)\leq
C(\xi).$$ In particular for $\xi_0={\rm diag}(1,0,0,...)$,
\begin{equation}\label{eq1}\int_{\Omega}(1-\varphi_{\omega}(\xi_0))\;\frac{m_t}{t}(d\omega)\leq C(\xi_0)=M.\end{equation} Let us recall that $\displaystyle\vp(\xi_0)=\Pi(\omega,1)=e^{-\gamma}\prod\limits_{k=1}^{\infty}{\frac
  {1}{1+\alpha_k}}.$ Let $\kappa_t$ denote the positive and bounded measure defined on $\Omega$ by $$\kappa_t=(1-\vp(\xi_0))\;\frac{m_t}{t}.$$
Since, for all $t>0$, $\kappa_t(\Omega)\leq M$, the set $\left\{\kappa_t\;|\;t>0\right\},$ is relatively compact for the weak topology $\sigma({\mathcal M}(\Omega),\mathscr
C_0(\Omega))$, where ${\mathcal M}(\Omega)$ is the set of positive and bounded measures on $\Omega$ and $\mathscr C_0(\Omega)$ is the set of continuous functions on $\Omega$ vanishing at $+\infty$. It follows that there exists a sequence $(t_j)$ in $]0, +\infty[$ converging to $0$, such that the measures $\kappa_{t_j}$ weakly converge to a positive and bounded measure $\kappa$, i.e. for all $f\in {\mathscr C}_0(\Omega)$,
$$\lim_{j\rightarrow\infty}\int_{\Omega}f(\omega)\,\kappa_{t_j}(d\omega)=
\int_{\Omega}f(\omega)\,\kappa(d\omega).$$
On the other hand, we have
\begin{equation}\label{eq.}\frac{1-e^{-t_j\psi(\xi)}}{t_j}=\int_{\Omega}\left[\frac{1-\vp(\xi)}{1-\vp(\xi_0)}-1\right]\,\kappa_{t_j}(d\omega)  + \frac{1-e^{-t_j\psi(\xi_0)}}{t_j}.\end{equation}
Let us verify that, for $\xi\neq 0$, the function $f$ definied on $\Omega$ by $$f(\omega)=\left\{\begin{aligned}&\frac{1-\vp(\xi)}{1-\vp(\xi_0)}-1\quad\rm{si}\;\omega\neq
  0,\\
&{\rm tr}(\xi^2)-1\quad\quad\quad\quad\rm{si}\;\omega=0,
\end{aligned}\right.$$
 belongs to ${\mathscr C}_0(\Omega)$.

First, remark that $\vp(\xi_0)=1$ if and only if $\omega=0$. Therefore, the function $f$ is well defined and continuous on $\Omega\backslash\{0\}$. Moreover, by Lemma \ref{lemm},
$$\lim_{\omega\to 0}\frac{1-\vp(\xi)}{1-\vp(\xi_0)}=\tr(\xi^2).$$ So, the function $f$ is continuous at $0$.
 We also have    $$\lim\limits_{\omega\to+\infty}\vp(\xi)=0 \quad {\rm and\  so}\quad \lim_{\omega\to\infty}f(\omega)=0.$$ As $j$ goes to $+\infty$ in the equation (\ref{eq.}), we get for all $\xi\neq 0$,
$$\begin{aligned}\psi(\xi)&=\int_{\Omega}f(\omega)\,\kappa(d\omega)+\psi(\xi_0)\\&=(\tr(\xi^2)-1)\kappa(\{0\})+\int_{\Omega\backslash\{0\}}\left(\frac{1-\vp(\xi)}{1-\vp(\xi_0)}-1\right)\kappa(d\omega)+\psi(\xi_0).\end{aligned}$$
As $\xi$ goes to $0$, we get $\displaystyle\kappa(\Omega)=\psi(\xi_0)$ and so 
$$\psi(\xi)=\tr(\xi^2)\kappa(\{0\})+\int_{\Omega\backslash\{0\}}\frac{1-\vp(\xi)}{1-\vp(\xi_0)}\,\kappa(d\omega).$$
Let us write
$$\psi(\xi)=A\;\tr(\xi^2)+\int_{\Omega\backslash\{0\}}(1-\vp(\xi))\,\nu(d\omega)$$
with
 $A=\kappa(\{0\})$ and $\nu$ is the measure defined on $\Omega$ by $$\nu=\frac{1}{1-\vp(\xi_0)}\kappa_{\mid_{{\Omega\backslash\{0\}}}}.$$
Now, let us check that  $$\int_{{\Omega\backslash\{0\}}}\frac{||\omega||}{1+||\omega||}\,\nu(d\omega)<+\infty.$$

The function $\displaystyle\omega\mapsto\frac{||\omega||}{(1-\vp(\xi_0))(1+||\omega||)}$ is continuous on $\Omega\backslash\{0\}$. It has as limit $1$ at $0$ and $+\infty$. It is, in consequence, a bounded function and so
$$\int_{\Omega}\frac{||\omega||}{1+||\omega||}\nu(d\omega)=\int_{\Omega}\frac{||\omega||}{(1-\vp(\xi_0))(1+||\omega||)}\kappa(d\omega)<+\infty.$$


Finally, we obtain $$\psi(\xi)=A\tr(\xi^2)+\int_{\Omega\backslash\{0\}}(1-\vp(\xi))\,\nu(d\omega).$$ To get the uniqueness of $A$, we use the fact that the modified P{\'o}lya function $\Pi(\omega,.)$ is of positive type on $\Bbb R$, of class $\mathscr C^2$ and verifies $\Pi(\omega,0)=1$. Therefore
$$\begin{aligned}1-\Pi(\omega,s)=\int_\Bbb
  R(1-\cos(su))\mu(du)\leq\frac{1}{2}s^2\int_\Bbb Ru^2\mu(du)=-s^2\Pi''(\omega,0).
\end{aligned}$$ Besides, $$\vp(s\xi_0)=\Pi(\omega,s)\;\;\mbox{and}\;\; \Pi''(\omega,0)=-2||\omega||.$$
It follows that for all $s\geq 1$,
$$\frac{1-\vp(s\xi_0)}{s^2}\leq\left\{\begin{aligned}&2||\omega||\qquad\rm{si}\;
    ||\omega||\leq 1,\\&2\qquad\qquad\;\;\rm{si}\;||\omega||\geq 1.\end{aligned}\right.$$
By writing
$$\frac{\psi(s\xi_0)}{s^2}=A+\int_{\Omega\backslash\{0\}}\frac{1-\vp(s\xi_0)}{s^2}\,\nu(d\omega),$$
the dominated convergence theorem implies the uniqueness of the constant $A$ :
$$\lim_{s\to +\infty}\frac{\psi(s\xi_0)}{s^2}=A \ .$$

${\rm\bf c)}$ \ {\bf Uniqueness of the measure $\nu$} : Let $\nu_1$ and $\nu_2$ be two positive measures on $\Omega\backslash\{0\}$ such that, for all \;$\xi\in V_\infty$,

 \bq\int_{\Omega\backslash\{0\}}(1-\vp(\xi))\,\nu_1(d\omega)=\int_{\Omega\backslash\{0\}}(1-\vp(\xi))\,\nu_2(d\omega),\eq and \bq\int_{\Omega\backslash\{0\}}\frac{||\omega||}{1+||\omega||}\;\nu_i(d\omega)<\infty,\qquad(i=1,2).\vspace{.3cm}\eq  Then, by Proposition \ref{unique}, for all \;$\xi,\eta\in V_\infty$,

$$\int_{\Omega\backslash\{0\}}\vp(\eta)(1-\vp(\xi))\,\nu_1(d\omega)=\int_{\Omega\backslash\{0\}}\vp(\eta)(1-\vp(\xi))\,\nu_2(d\omega).$$ Let us put for all $\eta\in V_\infty$ $$\widehat\varphi_i(\eta)=\int_{\Omega\backslash\{0\}}\vp(\eta)(1-\vp(\xi))\,\nu_i(d\omega) \qquad (i=1, 2).$$ The function $\widehat\varphi_i$ is a continuous, $K_{\infty}$-invariant and positive definite function which is dominated by a $\nu_1$-int{e}grable function on $\Omega\backslash\{0\}$ :
 $$|\vp(\eta)(1-\vp(\xi))|\leq 1-\vp(\xi).$$

Let us put $$\widehat{\nu_{i,\xi}}(d\omega)=(1-\vp(\xi))\,\indic_{{}_{\Omega\backslash\{0\}}}(\omega)\,\nu_i(d\omega) \qquad (i=1, 2),$$
where $\indic_{{}_{\Omega\backslash\{0\}}}$ is the characteristic function of the set ${\Omega\backslash\{0\}}$. Then, for all $\eta\in
V_\infty$,$$\int_{\Omega}\vp(\eta)\,\widehat{\nu_{1,\xi}}(d\omega)=\int_{\Omega}\vp(\eta)\,\widehat{\nu_{2,\xi}}(d\omega).$$

By the uniqueness of the representing measure in the generalized Bochner theorem (Theorem \ref{th3}), we get $$\nu_1=\nu_2\qquad\mbox{on}\;\;\Omega\backslash\{0\}. \qquad \qquad \hfill\square$$}

\end{paragraph}





\begin{paragraph}{\bfseries{Acknowledgment}}{\rm The author is indebted to the referee for his valuable comments, remarks and suggestions which improved the presentation of this work.}

\end{paragraph}

\bibliographystyle{amsalpha}
 \vspace{0.5cm}

\begin{center}

\address{King Faisal University\hfill Institut préparatoire aux \ \ \ \ \ \ \ \ \ \ \ \ \\
Department of Mathematics\hfill \'Etudes d'Ingénieurs de Nabeul  \ \ \  \ \\ 
P.O.Box : 400, Al Hassa 31982 \hfill Campus Universitaire Merezka \ \ \ \ \\ 
Saudi Arabia \hfill 8000 Nabeul, Tunisia} \ \ \ \ \ \ \ \ \ \ \ \ \ \ \ \ \ \\

\email{mrabaoui@kfu.edu.sa} \hfill \email{Marouane.Rabaoui@ipein.rnu.tn} \  \ \\


\end{center}



\end{document}